\documentclass[12pt]{article}

\usepackage{latexsym}
\usepackage{bezier}

\begin{document}

%%%%%%%%%%%%%%%%%%%%%%%%%%%%%%%%%%%%%%%%%%%%%%%%%%%%%%%%%%%% begin lama.sty
%  LAMA.STY  (auszug)                                           6. M\"arz 92
%

% Spezielle mathematische Symbole ----------------------------------------
\newcommand{\EE}{\mathop{\rm I\! E}\nolimits}
\newcommand{\E}{\mathop{\rm E}\nolimits}
\newcommand{\I}{\mathop{\rm Im\, }\nolimits}
\newcommand{\Str}{\mathop{\rm Str\, }\nolimits}
\newcommand{\Sdet}{\mathop{\rm Sdet\, }\nolimits}
\newcommand{\STr}{\mathop{\rm STr\, }\nolimits}
\newcommand{\R}{\mathop{\rm Re\, }\nolimits}
\newcommand{\CC}{\mathop{\rm C\!\!\! I}\nolimits}
\newcommand{\FF}{\mathop{\rm I\! F}\nolimits}
\newcommand{\KK}{\mathop{\rm I\! K}\nolimits}
\newcommand{\LL}{\mathop{\rm I\! L}\nolimits}
\newcommand{\MM}{\mathop{\rm I\! M}\nolimits}
\newcommand{\NN}{\mathop{\rm I\! N}\nolimits}
\newcommand{\PP}{\mathop{\rm I\! P}\nolimits}
\newcommand{\QQ}{\mathop{\rm I\! Q}\nolimits}
\newcommand{\RR}{\mathop{\rm I\! R}\nolimits}
\newcommand{\ZZ}{\mathop{\sl Z\!\!Z}\nolimits}
% ------------------------------------------------------------------------
\newcommand{\integer}{\mathop{\rm int}\nolimits}
\newcommand{\erf}{\mathop{\rm erf}\nolimits}
\newcommand{\diag}{\mathop{\rm diag}\nolimits}
\newcommand{\fl}{\mathop{\rm fl}\nolimits}
\newcommand{\eps}{\mathop{\rm eps}\nolimits}
\newcommand{\var}{\mathop{\rm var}\nolimits}

%%%%%%%%%%%%%%%%%%%%%%%%%%%%%%%%%%%%%%%%%%%%%%%%%%%%%%%%%%%% end lama.sty

\newcommand{\pfeil}{\rightarrow}

\newcommand{\kat}{{\cal C}}
\newcommand{\rmat}{{\cal R}}
\newcommand{\oalg}{{\cal A}}
\newcommand{\falg}{{\cal F}}
\newcommand{\eich}{{\cal G}}
\newcommand{\hilb}{{\cal H}}
\newcommand{\calm}{{\cal M}}
\newcommand{\mod}{{\cal M}}
\newcommand{\kegel}{{\cal O}}
\newcommand{\kegels}{{\cal K}}
\newcommand{\bigrho}{\rho_\oplus}
\newcommand{\bigphi}{\phi_\oplus}

\newcommand{\tprod}{\otimes}

\newcommand{\horab}{\rule[-1mm]{0pt}{5mm}}
\newcommand{\iso}{\stackrel{\sim}{=}}
\newcommand{\quer}[1]{\overline{#1}}
\newcommand{\schlange}[1]{\widetilde{#1}}
\newcommand{\CVO}[3]{{{#3 \choose #1\hspace{3pt}#2}}}
\newcommand{\clebsch}[6]{{\left[{#1\atop#4}\hspace{3pt}{#2\atop#5}
                                     \hspace{3pt}{#3\atop#6}\right]}}
\newcommand{\sjsymbol}[6]{{\left\{{#1\atop#4}\hspace{3pt}{#2\atop#5}
                                     \hspace{3pt}{#3\atop#6}\right\}}}
\newcommand{\spann}{{{\rm span}}}
\newcommand{\Nat}{{{\rm Nat}}}
\newcommand{\Mor}{{{\rm Mor}}}
\newcommand{\End}{{{\rm End}}}
\newcommand{\ev}{{{\rm ev}}}
\newcommand{\coev}{{{\rm coev}}}
\newcommand{\id}{{{\rm id}}}
\newcommand{\Id}{{{\rm Id}}}
\newcommand{\Vec}{{{\rm Vec}}}
\newcommand{\Rep}{{{\rm Rep}}}
\newcommand{\ket}[1]{|#1\rangle}
\newcommand{\vak}{\ket{0}}
\newcommand{\bild}[3]{{        % \bild{H"ohe}{Unterschrift}{ref-label}
  \unitlength1mm
  \begin{figure}[ht]
  \begin{picture}(120,#1)\end{picture}
  \caption{\label{#3}#2}
  \end{figure}
}}

\newcommand{\zopf}{{\begin{picture}(15,10)
                    \qbezier(5,10)(5,7.5)(10,5)
                    \qbezier(10,5)(15,2.5)(15,0)
                    \qbezier(11,6)(15,8)(15,10)
                    \qbezier(5,0)(5,2)(9,4)
                    \end{picture}}}
\newcommand{\invzopf}{{\begin{picture}(15,10)
                    \qbezier(5,10)(5,8)(9,6)
                    \qbezier(11,4)(15,2)(15,0)
                    \qbezier(10,5)(15,7.5)(15,10)
                    \qbezier(5,0)(5,2.5)(10,5)
                    \end{picture}}}
\newcommand{\symzopf}{{\begin{picture}(15,10)
                    \qbezier(5,10)(5,7.5)(10,5)
                    \qbezier(10,5)(15,2.5)(15,0)
                    \qbezier(10,5)(15,7.5)(15,10)
                    \qbezier(5,0)(5,2.5)(10,5)
                    \end{picture}}}
\newcommand{\evpic}{{\begin{picture}(15,10)
                   \put(10,5){\oval(10,10)[b]}
                   \put(5,5){\line(0,1){5}}
                   \put(15,5){\line(0,1){5}}
                     \end{picture}}}
\newcommand{\coevpic}{{\begin{picture}(15,10)
                   \put(10,5){\oval(10,10)[t]}
                   \put(5,5){\line(0,-1){5}}
                   \put(15,5){\line(0,-1){5}}
                     \end{picture}}}
\newcommand{\evbox}[1]{{\begin{picture}(16,10)
                   \put(5,5){\line(0,1){5}}
                   \put(15,5){\line(0,1){5}}
                   \put(4,0){\framebox(12,5){{\small #1}}}
                         \end{picture}}}
\newcommand{\coevbox}[1]{{\begin{picture}(16,10)
                   \put(5,5){\line(0,-1){5}}
                   \put(15,5){\line(0,-1){5}}
                   \put(4,5){\framebox(12,5){{\small #1}}}
                         \end{picture}}}
\newcommand{\morbox}[1]{{\begin{picture}(10,10)
                   \put(5,2){\line(0,-1){2}}
                   \put(5,8){\line(0,1){2}}
                   \put(1.5,2){\framebox(7,6){{\small #1}}}
                         \end{picture}}}
\newcommand{\bigmorbox}[1]{{\begin{picture}(10,10)
                   \put(5,2){\line(0,-1){2}}
                   \put(5,8){\line(0,1){2}}
                   \put(0.5,2){\framebox(9,6){{\small #1}}}
                         \end{picture}}}
\newcommand{\hugemorbox}[1]{{\begin{picture}(15,10)
                   \put(7.5,2){\line(0,-1){2}}
                   \put(7.5,8){\line(0,1){2}}
                   \put(1,2){\framebox(13,6){{\small #1}}}
                         \end{picture}}}
\newcommand{\vertexbox}[1]{{\begin{picture}(16,10)
                   \put(15,8){\line(0,1){2}}
                   \put(5,8){\line(0,1){2}}
                   \put(10,2){\line(0,-1){2}}
                   \put(4,2){\framebox(12,6){{\small #1}}}
                         \end{picture}}}
\newcommand{\invvertexbox}[1]{{\begin{picture}(16,10)
                   \put(15,2){\line(0,-1){2}}
                   \put(5,2){\line(0,-1){2}}
                   \put(10,8){\line(0,1){2}}
                   \put(4,2){\framebox(12,6){{\small #1}}}
                         \end{picture}}}

\newenvironment{bew}{Proof:}{\hfill$\Box$}

\newtheorem{bem}{Remark}
\newtheorem{bsp}{Example}
\newtheorem{axiom}{Axiom}
\newtheorem{de}{Definition}
\newtheorem{satz}{Proposition}
\newtheorem{lemma}[satz]{Lemma}
\newtheorem{kor}[satz]{Corollary}
\newtheorem{theo}[satz]{Theorem}

\newcommand{\sbegin}[1]{\small\begin{#1}}
\newcommand{\send}[1]{\end{#1}\normalsize}

\sloppy
%%%%%%%%%%%%%%%%%%%%%%%%%%%%%%%%%%%%%%%%%%%%%%%%%%%%%%%%%%%%%%%%%%%%%%%%%

\title{Reconstruction of Weak Quasi Hopf Algebras}
\author{Reinhard H\"aring\\
Mathematisches Institut, SFB 170\\ Bunsenstr. 3-5\\
 37073 G\"ottingen, Germany\\
email: haering@cfgauss.uni-math.gwdg.de}
\date{March, 14, 1995}
\maketitle

\begin{abstract}

All rational semisimple braided tensor categories are representation
categories of weak quasi Hopf algebras. To proof this result we
construct for any given category of this kind a
weak quasi tensor functor to the category of finite dimensional
vector spaces. This allows to reconstruct  a weak quasi Hopf algebra
with the given category as its representation category.
\end{abstract}

{\small \tableofcontents}

\section{Introduction}

Semisimple braided tensor categories are the structure underlieing
the quantum invariants of links and 3-manifolds \cite{turaev}.
The most useful examples are derived from the (nonsemisimple)
representation categories of quantum groups by elimination
of not fully decomposable objects and nilpotent morphisms.

These cleaned up versions are no longer  representation
categories of usual quantum groups.
It is the purpose of this paper to show that they nevertheless
arise as representation categories of appropriate algebras.
It is always possible to reconstruct a weak quasi Hopf algebras,
as introduced by Mack and Schomerus \cite{ms1}, that has the given
category as its representation category.

This result is established in two steps. First we define the notion
of a weak quasi tensor functor and show by construction that for any
rational braided semisimple tensor category $\kat$ such a functor $F$
to the category of finite dimensional vector spaces exists.

If this functor was a tensor functor in the usual sense then
Majid's reconstruction theorem \cite{majid1,majid2} would be
applicable and would assert the existence of an
associated  quasi Hopf algebra. It is fairly easy to show
\cite{kratz,nill} that such tensor functors don't exist for a
large classes of rational semisimple braided tensor categories.

However Majid's lines of thought can be applied even to the case
of a weak quasi tensor functor $F:\kat\pfeil\Vec$. This
generalized reconstruction theorem (section \ref{recth})
shows that the set $\Nat(F,F)$ of natural transformations
from $F$ to itself can be equipped with the structure of a weak
quasi Hopf algebra $H=H(\kat,F)$ such that $F$ factors over
$\Rep(H)$, i.e. there is a tensor functor $G:\kat\pfeil\Rep(H)$
such that $F=V\circ G$ where $V:\Rep(H)\pfeil\Vec$ is the
forgetful functor which assigns to any representation its
underlieing vector space.

Combining this reconstruction theorem with the construction
of weak quasi tensor functors we conclude that
every rational semisimple rigid braided tensor category
is the representation category of some weak quasi Hopf algebra.

\section{Braided Tensor Categories}

\subsection{Definitions}

The objects of a category $\kat$ are denoted by $X\in{\rm Obj}(\kat)$,
the morphisms  between
$X,Y\in {\rm Obj}(\kat)$ with ${\rm Mor}(X,Y)$.
We use the shorthand ${\rm End}(X):={\rm Mor}(X,X)$.
The identity functor of a category well be denoted by $\Id$ and
the set of natural transformations between two functors by
$\Nat(F,G)$.

A category $\kat$ is called {\em\bf monoidal}  if
there is a functor $\tprod:\kat\times\kat\pfeil\kat$ together
with a functorial isomorphism
$\Phi_{X,Y,Z}:X\tprod(Y\tprod Z)\stackrel{\sim}{\longrightarrow}
(X\tprod Y)\tprod Z$ satisfying the
{\em\bf pentagon identity}:
$\Phi\Phi=(\Phi\tprod 1)\Phi(1\tprod \Phi)$ and
an {\em\bf identity object} $1\in {\rm Obj}(\kat)$, such that
$r_X:X\mapsto 1\tprod X$ and $l_X:X\mapsto X\tprod1$ are
equivalences of categories compatible with $\Phi$:
$\Phi_{1,X,Y}\circ l_{X\tprod Y}=l_X\tprod \id_Y$.
It is called {\em\bf strict} if $1\tprod X=X\tprod 1=X,
(X\tprod Y)\tprod Z=X\tprod(Y\tprod Z),\Phi=\id,r=l=\id$.
Thanks to MacLanes coherence theorem all equivalence classes of
monoidal categories include strict ones.
A monoidal category is called {\em\bf braided tensor category}
if there is a functorial isomorphism
$\Psi_{X,Y}:X\tprod Y\stackrel{\sim}{\longrightarrow}Y\tprod X
\; X,Y\in {\rm Obj}(\kat) $
satisfying the two
{\em\bf hexagon identities}
$\Phi\Psi\Phi=(\Psi\tprod1)\Phi(1\tprod\Psi)$ and
$\Phi^{-1}\Psi\Phi^{-1}=(1\tprod\Psi)\Phi(\Psi\tprod1)$ as well as
$\Phi(r\tprod1)=(1\tprod\Psi)(1\tprod r)$ and $l_X=\Psi\circ r_X$.
In a {\em\bf tensor category}
or (in contrast to braided  tensor categories)
{\em\bf symmetric tensor category}
the identity $\Psi_{X,Y}\Psi_{Y,X}=\id_{X\tprod Y}$ should hold.

We assume all categories to be abelian (and all functors to be additive)
with direct sum $\oplus$ and zero element $0$. Then $\End(1)$ is a ring
and we assume it in addition to be a field which we denote by $\KK$.

An object $X\in{\rm Obj}(\kat)$ is called
{\em\bf indecomposable} if
${\rm End}(X)=\spann\;\id_X\oplus {\cal N}$ where
${\cal N}$ consists only of nilpotent elements, and
$X$ is called {\em\bf irreducible},
if ${\cal N}=0$. The set of irreducible objects is denoted by
{\bf ${\rm Obj}_{irr}$}.
In a {\em\bf fully reducible} category
all $X\in{\rm Obj}(\kat)$ are isomorphic to sums of irreducible objects.

Let  $\nabla\subset{\rm Obj}(\kat)$ denote a set containing
one object out of every equivalence class of irreducible objects.

In a {\em\bf quasi rational category} every object
is isomorphic to a finite sum of indecomposable objects.
A {\em\bf rational category} is a quasi rational category
with only finitely many equivalence classes
of indecomposable objects.
$\kat$ is called {\em\bf irredundant},
if $X\iso Y\Rightarrow X=Y$ and it is called
{\em\bf locally finite} if all $\Mor(X,Y)$ are finite dimensional
vector spaces (Note that quasi rational
categories are locally finite.).
A locally finite abelian  braided tensor category
is called {\em\bf semisimple},
if all ${\rm End}(X)$ are semisimple algebras.
By Wedderburn's theorem we get for locally finite categories
the equivalence of semisimplicity and full reducibility.

In a {\em\bf $C^\ast$-category} $\kat$
all ${\rm Mor}(X,Y)$ are Banach spaces with an
antilinear involution $\dagger:{\rm Mor}(X,Y)\pfeil {\rm Mor}(Y,X)$
such that $(fg)^\dagger=g^\dagger f^\dagger,
||f^\dagger f||=||f||^2$ (This implies that ${\rm End}(X)$
is a unital $C^\ast$-algebra.) and $\Psi^\dagger=\Psi^{-1},
\Phi^\dagger=\Phi^{-1}$.

A functor $F:\kat_1\pfeil \kat_2$ is called {\rm\bf faithful} if
$F:{\rm Mor}(X,Y)\pfeil {\rm Mor}(F(X),F(Y)) $
is injective for all $X,Y\in{\rm Obj}(\kat_1)$.

A tensor category is {\em\bf rigid} if it has dual objects, i.e.
there is a map of objects
$X\mapsto X^\ast$ and morphisms $\ev_X\in{\rm Mor}(X^\ast\tprod X,1)$
({\em\bf evaluation})
and $\coev_X\in{\rm Mor}(1,X\tprod X^\ast)$ ({\em\bf coevaluation}) with the
properties
\begin{equation}\label{evcoevax}
(\id\tprod \ev_X)\circ(\coev_X\tprod \id)=\id_X\qquad
(\ev_X\tprod \id)\circ(\id\tprod \coev_X)=\id_{X^\ast}
\end{equation}
The object mapping $\ast$ extends naturally
to an involutive (in the sense that $\ast\circ\ast$
is equivalent to $Id$.) contravariant functor by the definition
$f^\ast:=(\ev_Y\tprod \id)\circ(\id\tprod f\tprod \id)\circ(\id\tprod \coev_X)
\in{\rm Mor}(Y^\ast,X^\ast)$ for $f\in{\rm Mor}(X,Y)$.

Evaluation and coevaluation are unique up to a unique isomorphism.
It is always assumed that $\coev$ is an isometry if $\kat$ is a
$C^\ast$-category.
The duality map $\ast$ induces an involution $X\mapsto\hat{X}$ of $\nabla$
such that $X^\ast\iso \hat{X}$ for all $X\in\nabla$.

A {\em\bf ribbon category} is a braided category with the additional
structure of a twist, i.e. natural isomorphisms
$\sigma\in{\rm Nat}(\Id,\Id)$ obeying
$\sigma(X)^\ast=\sigma(X^\ast)$ and
$\Psi_{Y,X}\circ\Psi_{X,Y}\circ\sigma(X\tprod Y)=\sigma(X)\tprod\sigma(Y)$.
Note that these axioms imply
$(\sigma(X)^2\tprod \id)\circ \coev_X=
 \Psi_{X^\ast,X}\circ\Psi_{X,X^\ast}\circ \coev_X$
which shows that in symmetric tensor categories one has
$\sigma(X)^2=\id_X$.

Every $C^\ast$ category carries a natural ribbon structure: Define the
{\em\bf statistical parameter} by
$\lambda(X):= (\id_X\tprod \coev_X^\dagger)\circ
(\id_X\tprod\Psi_{X,X^\ast}^\dagger)\circ(\coev_X\tprod \id_X)\in{\rm End}(X)$
and write its polar decomposition as
$\lambda(X)=\sigma(X)^{-1}\circ P(X)$ where
the positive part is $P(X)$ and the unitary part $\sigma(X)^{-1}$
yields the desired ribbon structure.

A ribbon structure allows the definition of a trace on ${\rm End}(X)$ by
$tr_X(f):=\ev_X\circ(\id\tprod(f\circ\sigma(X)^{-1}))
\circ\Psi_{X,X^\ast}\circ \coev_X$ and a dimension
$d(X):= tr(\id_X)$. One has
$tr(f\circ g)=tr(g\circ f), tr(f)=tr(f^\ast),tr(f\tprod g)=tr(f)tr(g),
tr(f\oplus g)=tr(f)+tr(g), tr(f^\dagger)=tr(f)^\dagger$ and hence
$d(X\oplus Y)=d(X)+d(Y), d(X\tprod Y)=d(X)d(Y), d(X^\ast)=d(X)$.
Related to the trace is the {\em\bf conditional expectation}
$E_X:{\rm End}(A\tprod X)\pfeil{\rm End}(A), f\mapsto
(\id_A\tprod (\ev_X\circ\Psi_{X,X^\ast}))\circ(f\tprod \id_{X^\ast})\circ
(\id_A\tprod \sigma(X)^{-1}\tprod \id_{X^\ast})\circ(\id\tprod \coev_X)$.
Easy calculations show that it obeys
$E_X((g\tprod \id)\circ h\circ(f\tprod \id))=g\circ E_X(h)\circ f,
E_X(f\tprod \id_X)=d(X)f,E_X(f\tprod\Psi_{X,X})=f\tprod\sigma(X)^{-1}$.
On ${\rm End}(X^{\tprod n})$ the iterated expectation $E_X^n$
coincides with $tr$.

In the application we have in mind, the categories are
representation categories of algebras.

For an algebra $A$ we let $\Rep (A)$ denote
its {\em\bf representation category}.
The objects are the  representations of $A$
(it is common to consider only  a special class of representations)
and the morphisms are the intertwiners.

$\Rep (A)$ is a braided tensor category  if $A$ permits products
of representations which are symmetric up to isomorphisms.

Examples of  monoidal categories are superselection categories in
algebraic quantum field theory, the Moore/Seiberg categories in conformal
quantum field theory
and the representation categories of quantum groups.
For any quasitriangular Hopf algebra $H$ the class of its
representations is a braided tensor
category with braid isomorphism between two representations
$\varrho_1,\varrho_2$ naturally given by
\begin{equation}\label{zopfhopf}
\Psi_{\varrho_1,\varrho_2}(v_1\tprod v_2):=\tau\circ
  (\varrho_1\otimes\varrho_2)(R)(v_1\tprod v_2)
\end{equation}
Here $\tau$ is the flip operator $a\otimes b\mapsto b\otimes a$ and
$R$ is the usual $R$-matrix, i.e. $\tau\circ\Delta(a)=R\Delta(a)R^{-1}$.

The subcategory $\Rep(H)^{fd}$
of finite dimensional representations is rigid
thanks to the conjugate representation. Usually we will consider only this
subcategory and hence omit the superscript $fd$.

\subsection{Description of semisimple categories via
polynomial equations}\label{peq}

Let $\kat$ denote a semisimple braided
tensor category. For each triple $X,Y,Z\in\nabla$ let
$N_{X,Y}^Z$ denote the dimension
of ${\rm Mor}(X\tprod Y,Z)$ and choose a basis
$\phi(e)\in{\rm Mor}(X\tprod Y,Z)$
($e=^i\CVO{X}{Y}{Z}$ is a multi index with
$i\in\{1,\ldots,N_{X,Y}^Z\}$.).
The composite morphisms
$\phi^i\CVO{X}{Y}{Z}\circ\Psi_{Y,X}\in{\rm Mor}(Y\tprod X,Z)$ and
$\phi^i\CVO{X}{M}{R}\circ(\id_X\tprod\phi^j\CVO{Y}{Z}{M})
\in{\rm Mor}(X\tprod(Y\tprod Z))$
can then be expanded in the basis via matrices
\begin{eqnarray}
\phi(e)\circ\Psi&=&\sum_f\Omega_{e,f}\phi(f)\\
\phi(e_2)(\id\otimes\phi(e_1))&=&
\sum_{e,f}F_{e_1,e_2;f,e}\phi(e)(\phi(f)\otimes \id)
\end{eqnarray}
It follows straightforward from the axioms of braided tensor categories that
these matrices satisfy the Moore/Seiberg polynomial equations
\cite{moores2}.

Two semisimple rigid braided tensor categories are
equivalent if they are equivalent as ordinary
categories and they share the same structural data $\Omega,F$.

Moore/Seiberg have shown \cite{moores2} that in
the opposite direction every solution to their equations yields
such a category. Their construction is essentially the following:
Take a set of irreducible objects $X_i,i\in{\cal I}$ and set
${\rm Mor}(X_i,X_j):=\KK\delta_{i,j}\id_{X_i}$. Tensor products
are formally introduced via
$X_i\tprod X_j:=\bigoplus_l V_{i,j}^l\otimes X_l$ where $V_{i,j}^l$ are
$N_{i,j}^l$ dimensional vector spaces of morphisms ${\rm Mor}
(X_i\tprod X_j,X_l)$.
The braid isomorphism operates
on this tensor product via the operation of $\Omega$ on $V_{i,j}^l$.

\subsection{Graphical Calculus}

There is a handy notation for visualizing morphisms in strict braided
tensor categories.  Some basic morphisms are shown in figure
\ref{graphdef}).
The tensor product $f\tprod g$ is displayed by
drawing the picture for $f$ to the left of the picture for $g$
while $f\circ g$ is visualized by placing the picture for $g$ on top
of the picture for $f$. The unit object and its identity morphism
are usually not displayed.

\setlength{\unitlength}{1mm}

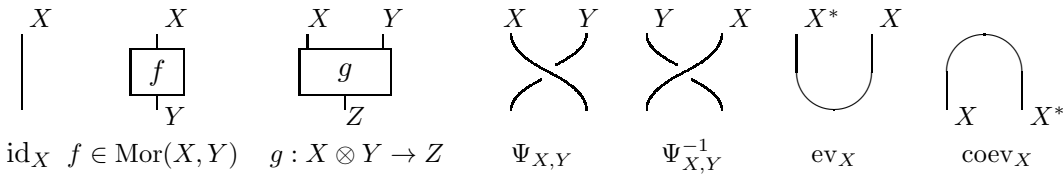
\begin{figure}[ht]
\begin{picture}(150,30)
\put(3,21){{\footnotesize $X$}}
\put(2,10){\line(0,1){10}}
\put(0,3){{\small $\id_X$}}

\put(21,21){{\footnotesize $X$}}
\put(21,8){{\footnotesize $Y$}}
\put(15,10){\morbox{$f$}}
\put(8,3){{\footnotesize $f\in\Mor(X,Y)$}}

\put(35,10){\vertexbox{$g$}}
\put(45,8){{\footnotesize $Z$}}
\put(40,21){{\footnotesize $X$}}
\put(50,21){{\footnotesize $Y$}}
\put(35,3){{\footnotesize $g:X\tprod Y\pfeil Z$}}

\put(62,10){\zopf}
\put(66,21){{\footnotesize $X$}}
\put(76,21){{\footnotesize $Y$}}
\put(67,3){{\footnotesize $\Psi_{X,Y}$}}
\put(80,10){\invzopf}
\put(86,21){{\footnotesize $Y$}}
\put(96,21){{\footnotesize $X$}}
\put(87,3){{\footnotesize $\Psi_{X,Y}^{-1}$}}

\put(100,10){\evpic}
\put(106,21){{\footnotesize $X^\ast$}}
\put(116,21){{\footnotesize $X$}}
\put(107,3){{\footnotesize $\ev_{X}$}}
\put(120,10){\coevpic}
\put(126,8){{\footnotesize $X$}}
\put(136,8){{\footnotesize $X^\ast$}}
\put(127,3){{\footnotesize $\coev_{X}$}}

\end{picture}
\caption{\label{graphdef}
Graphical notations for morphisms in braided tensor categories}
\end{figure}

%\begin{figure}[h]
%\begin{picture}(150,60)
%\end{picture}
%\caption{Axioms for braided tensor categories}
%\end{figure}

\subsection{Weak Tensor Functors}

A functor $F:\kat_1\pfeil\kat_2$ between two monoidal
categories is called
{\em\bf monoidal} (resp. {\em\bf weak monoidal})
if there is a functorial isomorphism (resp. epimorphism) $c_{X,Y}$
\begin{equation}
c_{X,Y}:F(X)\tprod_2F(Y)\stackrel{\sim}{\longrightarrow}F(X\tprod_1Y)
\end{equation}
such that $F$ becomes compatible with the associator and the unit:
\begin{equation}\label{assod}\begin{array}{ccccc}
F(X)\tprod(F(Y)\tprod F(Z))&\stackrel{1\tprod c}{\longrightarrow}&
F(X)\tprod F(Y\tprod Z)&\stackrel{c}{\longrightarrow}&F(X\tprod(Y\tprod Z))\\
\left\downarrow \horab\Phi_2\right. &&&&
\left\downarrow \horab F(\Phi_1)\right. \\
(F(X)\tprod F(Y))\tprod F(Z)&\stackrel{c\tprod1}{\longrightarrow}&
F(X\tprod Y)\tprod F(Z)&\stackrel{c}{\longrightarrow}&F((X\tprod Y)\tprod Z)
\end{array}\end{equation}
\begin{equation} l_2|_{F({\rm Obj}(\kat_1))}=
  c^{-1}\circ F(l_1):F(X)\mapsto F(1)\tprod_2 F(X)
  \iso1\tprod_2 F(X)                        \label{einskomp}
\end{equation}
A functor between two (braided) tensor categories is called
{\em\bf symmetric}  if it is
compatible with the braid isomorphism, i.e.
for all $X,Y\in{\rm Obj}(\kat)$ the diagram
\begin{equation}\begin{array}{ccc}
F(X)\tprod F(Y)&\stackrel{c}{\longrightarrow}&F( X\tprod Y)\nonumber\\
\left\downarrow\horab \Psi_2\right. &   &
\left\downarrow\horab F(\Psi_1)\right.\nonumber\\
F(Y)\tprod F(X)&\stackrel{c}{\longrightarrow}&
F(Y\tprod X)
\end{array}
\label{zopfkomp}
\end{equation}
is commutative.
A monoidal functor between   braided tensor categories is called
a {\em\bf tensor functor} if:
\begin{equation}
F(\Psi(X\tprod Y))\iso F(X\tprod Y) \label{zopfkomp2}
\end{equation}
(This property follows in all cases with exception
of the ultraweak case from the other axioms. One could therefore
formulate most of the present paper using only the term monoidal functor.)

If  (\ref{assod})  is not required $F$ is called
{\em\bf quasi tensor functor} and if
$c_{X,Y}$ is only an epimorphism (but with $c_{X,1}$ and $c_{1,X}$
remaining isomorphisms)
with right inverse $c_{X,Y}^{-1}$ then $F$ is only a
{\em\bf weak quasi tensor functor}.
Finally $F$ is called {\em\bf ultra weak quasi tensor functor}
if (\ref{einskomp})
and $c_{1,X},c_{X,1}$ are not postulated to be isomorphisms but
$c_{1,X}=c_{X,1}\circ\Psi^{\kat_2}$.

If $\kat_1$ and $\kat_2$ are rigid
then we demand in addition
the existence of functorial isomorphisms $d_X:F(X)^\ast\pfeil F(X^\ast)$.

If both categories are $C^\ast$ then $F$ is called {\em\bf isometric} if
$F(f^\dagger)=F(f)^\dagger$. Consistency with the tensor product requieres
then $c^{-1}$ to be an isometry (This is implied by the following calculation:
$c\circ(F(f^\dagger)\tprod F(g^\dagger))\circ c^{-1}=
F(f^\dagger\tprod g^\dagger)=
F((f\tprod g)^\dagger)=F(f\tprod g)^\dagger=
(c\circ(F(f)\tprod F(g))\circ c^{-1})^\dagger=
c^{-1\dagger}\circ(F(f^\dagger)\tprod F(g^\dagger))\circ c^\dagger
$).

$\kat_1$ and $\kat_2$ are {\em\bf equivalent as braided tensor categories} if
they are equivalent as usual categories with symmetric tensor
functors.

\subsection{Construction of (weak) Quasi Tensor Functors}
\label{ktab007}

\begin{de} A function defined on the irreducible objects of a
semisimple, rigid braided tensor category
$D:{\rm Obj}_{irr}(\kat)\pfeil\NN_0$ which is
constant on equivalence classes is called
{\em\bf weak dimension function}, if:
\begin{equation}
D(1)=1, D(X)=D(X^\ast), D(X)D(Y)\geq\sum_{Z\in\nabla}
D(Z){\rm dim}({\rm Mor}(X\tprod Y,Z))
\end{equation}
$D$ is called {\em\bf dimension function} if equality holds.
\end{de}

Dimension functions allow the construction of  monoidal functors:

\begin{satz}\label{funktorkon}
Let $\kat$ be a quasi-rational semisimple,
rigid, braided tensor category and $D:{\rm Obj}(\kat)\pfeil\NN$
a (weak) dimension function.
Then there is a faithful (weak) quasi tensor functor
$F:\kat\pfeil\Vec$ into the category of
finite dimensional vector spaces.
\end{satz}

The following lemma will be used in the proof of the proposition.
\begin{lemma} Let $X\in {\rm Obj}(\kat)$ be an irreducible object
in a semisimple category. Then for all
$Y\in{\rm Obj}(\kat)$ we have
${\rm Mor}(Y,X)\iso{\rm Mor}(X,Y)^{\ast_\Vec}$
\end{lemma}
\begin{bew}
Let $g\in{\rm Mor}(Y,X)$ and define $\lambda_g\in
{\rm Mor}(X,Y)^{\ast_\Vec}$ by
$\lambda_g(f):= g\circ f\in{\rm Mor}(X,X)$.
This pairing is nondegenerate: Assume $g\neq 0$. Then, by
semisimplicity, $Y\iso X\oplus Y_1$. But this implies
existence of a $f\in{\rm Mor}(Y,X)$ such that $g\circ f=\id\neq 0$.
\end{bew}\\
\begin{bew}
For $X\in\nabla$ let $F(X):=\KK^{D(X)}$
and for arbitrary objects $Y\in{\rm Obj}(\kat)$ this is extended via
$F(Y):=\bigoplus_{X\in\nabla}{\rm Mor}(X,Y)\otimes F(X)$.
$F$ has to map morphisms $f\in{\rm Mor}(Y_1,Y_2)$ to morphisms
$F(f)\in{\rm Mor}(F(Y_1),F(Y_2))$.
Because of linearity, $F(f)$ needs only be defined on the
summands of type
${\rm Mor}(X,Y_1)\otimes F(X)$.
Let $F(f)(g\otimes x):=f\circ g\otimes x,x\in F(X)$ for
$g\in{\rm Mor}(X,Y_1)$

Assume $f_1,f_2\in{\rm Mor}(Y_1,Y_2), F(f_1)=F(f_2)$.
By the definition of $F$ this implies that for all
$X\in\nabla$ and for all $g\in{\rm Mor}(X,Y_1)$ we have
$f_1\circ g=f_2\circ g$.
Since $\kat$ is assumed to be semisimple we have an isomorphism
$\phi\in{\rm Mor}(X_{i_1}\oplus\ldots\oplus X_{i_n},Y_1)$
with $X_{i_l}\in\nabla$.
{}From this we get $p_{i_l}\in
{\rm Mor}(X_{i_1}\oplus\ldots\oplus X_{i_n},X_{i_l}), q_{i_l}\in
{\rm Mor}(X_{i_l},Y_1)$ such that $\phi=\sum_l q_{i_l}\circ p_{i_l}$.
Now $\phi$ is epi and we have $f_1\circ q_{i_l}=f_2\circ q_{i_l}$
by the above remark. Hence $f_1\circ\phi=f_2\circ\phi$
and by this $f_1=f_2$: $F$ is faithful.

$F$ satisfies $F(Y^\ast)\iso F(Y)^\ast$:
\begin{eqnarray*}
 F(Y^\ast)=\bigoplus_{X\in\nabla}
 {\rm Mor}(X,Y^\ast)\otimes F(X)
          \iso\bigoplus_{X\in\nabla}{\rm Mor}
          (X^\ast,Y^\ast)\otimes F(X^\ast)\\
          \iso\bigoplus_{X\in\nabla}{\rm Mor}
          (X,Y)^\ast\otimes F(X^\ast)
       \iso\bigoplus_{X\in\nabla}{\rm Mor}
       (X,Y)^\ast\otimes F(X)^\ast =F(Y)^\ast
  \end{eqnarray*}
The lemma is used in the third step and
the fourth step uses the fact that
$F(X)$ and $F(X^\ast)$ are vector spaces of equal dimension.

For every pair of  irreducible objects $X_1,X_2\in\nabla$
we choose an arbitrary (epi/iso)morphism

\[C_{X_1,X_2}:F(X_1)\otimes F(X_2)\pfeil F(X_1\tprod X_2)=
\bigoplus_{X\in\nabla}{\rm Mor}(X,X_1\tprod X_2)\otimes F(X)\]

$c$ is defined as an extension of $C$:
\begin{eqnarray*}
c_{Y_1,Y_2}&:&F(Y_1)\otimes F(Y_2)\pfeil F(Y_1\tprod Y_2)\\
c_{Y_1,Y_2}&:&\left(\bigoplus_{X_1\in\nabla}{\rm Mor}(X_1,Y_1)
\otimes F(X_1)\right)
\otimes
\left(\bigoplus_{X_2\in\nabla}{\rm Mor}(X_2,Y_2)
\otimes F(X_2)\right)
\\&& \pfeil
\bigoplus_{X\in\nabla}{\rm Mor}(X,Y_1\tprod Y_2)\otimes F(X)\\
c_{Y_1,Y_2}&:=&\bigoplus_{X_1,X_2\in\nabla}
(\Gamma\otimes \id)\circ (\id\otimes \id\otimes C_{X_1,X_2})\circ\tau_{2,3}\\
\Gamma&:&{\rm Mor}(X_1,Y_1)\otimes{\rm Mor}(X_2,Y_2)
\otimes{\rm Mor}(X,X_1\tprod X_2)\pfeil {\rm Mor}(X,Y_1\tprod Y_2)\\
&&\Gamma(f_1\otimes f_2\otimes g):=(f_1\tprod f_2)\circ g
\end{eqnarray*}
$c$ behaves functorial, i.e. for
$f_i\in{\rm Mor}(Y_i,\schlange{Y_i}), i=1,2$ we have
$F(f_1\tprod f_2)\circ c_{Y_1,Y_2}=c_{\schlange{Y_1},\schlange{Y_2}}\circ
 (F(f_1)\otimes F(f_2))$.
To see this we introduce $v_i\in F(Y_i),i=1,2$ as
\[v_i=\bigoplus_{A_i\in\nabla} g^{(A_i)}\otimes x^{(A_i)}\qquad
x^{(A_i)}\in F(A_i),g^{(A_i)}\in{\rm Mor}(A_i,Y_i)\]
Using the definitions and the shorthand
$C_{A_1,A_2}(x^{(A_1)}\otimes x^{(A_2)})=\bigoplus_{B\in\nabla}
q^B_{A_1,A_2}\otimes x^B_{A_1,A_2} $
we get
\begin{eqnarray*}
\lefteqn{c_{\schlange{Y_1},\schlange{Y_2}}\circ
 (F(f_1)\otimes F(f_2))(v_1\otimes v_2)=}\\
& & =\bigoplus_{A_1,A_2\in\nabla}
 (\Gamma\otimes \id) f_1\circ g^{(A_1)}\otimes f_2\circ g^{(A_2)}\otimes
 C_{A_1,A_2}(x^{(A_1)}\otimes x^{(A_2)})\\
& & =\bigoplus_{A_1,A_2,B\in\nabla}
 (f_1\circ g^{(A_1)}\otimes f_2\circ g^{(A_2)})\circ
  q^B_{A_1,A_2}\otimes x^B_{A_1,A_2}\\
& & = F(f_1\tprod f_2)\circ\left(\bigoplus_{A_1,A_2,B\in\nabla}
 (g^{(A_1)}\otimes g^{(A_2)})\circ q^B_{A_1,A_2}\otimes x^B_{A_1,A_2}
 \right) \\
& & = F(f_1\tprod f_2)\circ c_{Y_1,Y_2}(v_1\otimes v_2)
\end{eqnarray*}
\end{bew}

% \label{kfunktor}\label{symalg}

\begin{bem} The functorial isomorphisms $d_Y:F(Y)^\ast\pfeil F(Y^\ast)$
can be displayed explicitly. First note that
\[F(Y^\ast)=\bigoplus_{X\in\nabla}{\rm Mor}(X,Y^\ast)\otimes F(X)
=\bigoplus_{X\in\nabla}{\rm Mor}(\hat{X},Y^\ast)\otimes F(\hat{X}) \]
We calculate $d_Y^{-1}$ operating on one summand
$u^\ast\otimes v\in{\rm Mor}(\hat{X},Y^\ast)\otimes F(\hat{X})
\subset F(Y^\ast)$.
We have to choose (Actually $z_X$ and $\schlange{d_X}$
are fixed by demanding $F(\ev_Y)\circ c_{Y^\ast,Y}\circ (d_Y\otimes \id)=
\ev_{F(Y)}$.) isomorphisms $z_X=\schlange{z_X}^\ast
\in{\rm Mor}(X^\ast,\hat{X})$ and
$\schlange{d}_X:F(X)^\ast\pfeil F(\hat{X})$
for all $X\in\nabla$.
Now we map
\[ u^\ast\otimes v\mapsto u^\ast\circ\schlange{z_X}^\ast\otimes v
=(\schlange{z_X}\circ u)^\ast\otimes v
\in\bigoplus{\rm Mor}(Y,X)^{\ast_\kat}\otimes F(\hat{X}) \]
With the techniques of the proceeding lemma this is
$\lambda_{\schlange{z_X}\circ u}\otimes v
\in {\rm Mor}(X,Y)^{\ast_\Vec}$. Finally applying
$\id\otimes\schlange{d_X}^{-1}$ yields
$\lambda_{\schlange{z_X}\circ u}\otimes\schlange{d_X}^{-1}(v)
\in\bigoplus{\rm Mor}(X,Y)^{\otimes_\Vec}\otimes F(X)^\ast
=F(Y)^\ast$

Using this description of $d_Y$ one can show
$d_W\circ F(f)^\ast=F(f^\ast)\circ d_Y$ for $f\in{\rm Mor}(W,Y)$.
\end{bem}

\begin{bem}
With arbitrary choices of the $C$ morphisms in the proof of the theorem
the constructed functor will in general not be compatible
with the associativity constraints in the sense of {\rm (\ref{assod})}.
For a strict (i.e. $\Phi=\id$) category {\rm (\ref{assod})} reads
\[c_{X\tprod Y,Z}\circ(c_{X,Y}\otimes \id)=
  c_{X,Y\tprod Y}\circ(\id\otimes c_{Y,Z}) \]
This equation can be interpreted as a non abelian two-cocycle
condition. 	We will take up this point later on.
\end{bem}

Proposition \ref{funktorkon} reduces the problem of finding a functor
to finding a dimension function.
This is possible:
\begin{satz}
On rational, semisimple, rigid, braided tensor categories
there exist always weak dimension functions.
\begin{equation}
D(1):=1\qquad D(X):={\rm dim}\bigoplus_{Y,Z\in\nabla}
{\rm Mor}(Y\tprod X,Z)=\sum_{i,j}N_{X,i}^j
\end{equation}
\end{satz}
Other possibilities are
$D(1):=1, D(X):=max_{I,J\neq0}\sum_K N^K_{I,J}$
and in the algebraic formulation of QFT \cite{Haag}
 $ D(\rho):={\rm dim(span}\{(\rho_I\rho,\rho_J)\mid
 \rho_I,\rho_J\in\nabla\})$

\begin{bew}
\begin{eqnarray*}
D(X)D(Y)=\left(\sum_{s,r}N_{X,s}^r\right)\left(\sum_{S,R}N_{Y,S}^R\right)
=\sum_{s,r,S,R}N_{X,s}^rN_{Y,S}^R\geq\\
\sum_{K,N,M}N_{X,N}^KN_{Y,K}^M=\sum_{K,N,M}N_{X,Y}^KN_{K,N}^M=D(X\tprod Y)
\end{eqnarray*}
\end{bew}

\subsection{Weak and ultra weak quasi Hopf algebras}

The structure of most rational semisimple tensor categories
does not allow non weak dimension functions \cite{kratz,nill}.
This results from the fact that ordinary quantum groups at roots
of unity have indecomposable representations of zero
(quantum) dimension $d$. They arise in the tensor product
decomposition of simple representations and spoil many
of the intented applications, e.g.
the interpretation of ordinary quantum groups as
gauge symmetry algebras is impossible.

To discard the indecomposable representations one has to allow
that the coproduct of unity, $\Delta(1)$, is not $1\otimes 1$,
but a projector on the fully decomposable part.
This is the idea of Mack/Schomerus
encoded in the definition of {\em\bf weak quasi
Hopf algebras} as modifications of Drinfeld's quasi Hopf algebras.
As those they are unital algebras $H$ together with a comultiplication
$\Delta:H\pfeil H\otimes H$, counit $\epsilon:H\pfeil\KK$ and
antipode $S:H\pfeil H$. The coproduct is commutative up to conjugation
by $R\in H\otimes H$ and associative up to conjugation by
$\phi\in H\otimes H\otimes H$, that is for all $h\in H$ one has
$\phi((\id\otimes\Delta)\circ\Delta(h))=
  ((\Delta\otimes\id)\circ\Delta(h))\phi$.

\begin{eqnarray}
\phi^{-1}\phi&=&(\id\otimes\Delta)\Delta(1)\label{ff1}\\
\phi\phi^{-1}&=&(\Delta\otimes \id)\Delta(1)\label{ff2}\\
RR^{-1}&=&\Delta'(1), \quad\Delta':=\tau\circ\Delta\label{ff3}\\
R^{-1}R&=&\Delta(1)\label{ff4}\\
(\id\otimes\id\otimes\epsilon)(\phi)&=&(\id\otimes\epsilon\otimes \id)(\phi)
=(\epsilon\otimes \id\otimes \id)(\phi)=\Delta(1)\label{ff5}
\end{eqnarray}
For the sake of completness we also recall Drinfeld's form of
the antipode axiom for quasi Hopf algebras. It states the existence
of two invertible elements $\alpha,\beta\in H$ such that the
following relations hold:
\begin{eqnarray}\label{antiax}
\epsilon(a)\alpha&=&\sum_i S(a_i^{(1)})\alpha a^{(2)}_i\quad\forall a\in H\\
\epsilon(a)\beta&=&\sum_i a^{(1)}_i\beta S(a_i^{(2)})\quad
  \forall a\in H\label{antiax2}\\
\mbox{with\quad} \Delta(a)&=&\sum_i a_i^{(1)}\otimes a_i^{(2)}\nonumber
\end{eqnarray}

Is there some kind of algebra generalizing the ((weak) quasi) quantum groups
and observable algebras of algebraic quantum field theory?
We believe that ultra weak quasi quantum groups
as introduced in \cite{haring} may provide an answer.

\begin{de}[Ultra weak quasi Hopf algebra]
Let $A$ denote an unital algebra.
An {\em\bf $A$ ultra weak quasi Hopf algabra} $H$
is a $A$ bialgebra $H$
(left and right multiplication are denoted by
$\mu_l:A\otimes H\pfeil H, \mu_r:H\otimes A\pfeil H$)
and algebra morphisms $\eta:A\pfeil H,\epsilon:H\pfeil A$ such that
all axioms of a weak quasi Hopf algebra are fulfilled with the exception of
unit/counit properties which are replaced by:
\begin{equation}
\mu_l(\epsilon\otimes \id)\Delta=\mu_r(\id\otimes \epsilon)\Delta=\id_H\qquad
m(\id\otimes\eta)=\mu_r\qquad m(\eta\otimes \id)=\mu_l
\end{equation}
\end{de}

\section{Reconstruction Theorems} \label{recth}

Historically the first reconstruction theorem was the famous
Tannaka-Krein theorem: Given a symmetric tensor category and a
faithful tensor functor to $\Vec$ there is a group
with the given category as representation category.
Majid proved reconstruction theorems for
quasitriangular Hopf algebras and quasi Hopf algebras.
A reconstruction theorem for weak quasi Hopf algebras was suggested by
Kerler without a proof.

The {\em\bf forgetful functor}
$V:\Rep(H)\pfeil\Vec$
assigns to each representation the underlieing vector space.

We start in lemma \ref{anfang} by reviewing Majid's reconstruction theorem
for quasi Hopf algebras. Starting point for his construction is the
set $\Nat(F,F)$ of natural transformations of $F$.
\begin{eqnarray*}
H&:=&H(\kat,F):=\Nat(F,F)=\\
&&\{ h:{\rm Obj}(\kat)\pfeil\End_\Vec \mid
 h_X\in\End(F(X)),\\
&&\qquad F(f)\circ h_X=h_Y\circ F(f)
\forall X,Y\in{\rm Obj}(\kat)\forall f\in\Mor(X,Y)\}
\end{eqnarray*}

\begin{lemma} \label{anfang}
$H$ is a quasitriangular (quasi) Hopf algebra
if $F$ is a (quasi) tensor functor.
\end{lemma}
\begin{bew}
$H$ is a vector space by pointwise addition.
The multiplication is also defined pointwise:
$(hg)_X:=h_X\circ g_X\quad X\in{\rm Obj}(\kat), h,g\in H$. The unit is
$X\mapsto 1_X=\id_{F(X)}$.
(The ultra weak case is handled in Lemma \ref{ultralemma}.)

In $\Vec$
the following relation holds:
${\rm End}(F(X))\otimes{\rm End}(F(Y))\iso
{\rm End}(F(X)\otimes F(Y))$ so that $H\otimes H$ is given by
functions in two variables $X,Y$ (i.e. we understand
the tensor product algebraically.),
which map to ${\rm End}(F(X)\otimes F(Y))$.
The coproduct $\Delta:H\pfeil H\otimes H$ is defined by:
\begin{equation}\Delta(h)_{X,Y}:=
c_{X,Y}^{-1}\circ h_{X\tprod Y}\circ c_{X,Y}\end{equation}
This is compatible with multiplication:
\begin{eqnarray*}
(\Delta(h)\Delta(g))_{X,Y}&=&\Delta(h)_{X,Y}\circ\Delta(g)_{X,Y}\\
&=&c^{-1}_{X,Y}\circ h_{X\tprod Y}\circ c_{X,Y}\circ
c_{X,Y}^{-1}\circ g_{X\tprod Y}\circ c_{X,Y}\\
&=&c^{-1}_{X,Y}\circ h_{X\tprod Y}\circ g_{X\tprod Y}\circ c_{X,Y}=
\Delta(hg)_{X,Y}
\end{eqnarray*}

The counit is $\epsilon:H\pfeil\KK,\epsilon(h):=h_1$.
\[((\id \otimes\epsilon)\Delta(h))_X=\Delta(h)_{X,1}=
c_{X,1}^{-1}\circ h_{X\tprod 1}\circ c_{X,1}=h_{X\tprod 1}=h_X\]

The associator
 $\phi\in H\otimes H\otimes H$ is given by
\begin{equation}
\phi_{X,Y,Z}:=(c_{X,Y}^{-1}\otimes\id)\circ c^{-1}_{X\tprod Y,Z}\circ
F(\Phi_{X,Y,Z})\circ c_{X,Y\tprod Z}\circ(\id\otimes c_{Y,Z})
\end{equation}
For tensor functors this is trivial because of (\ref{assod}).
For quasi tensor functors it is invertible.
\begin{eqnarray*}
\lefteqn{(\phi(1\otimes\Delta)\Delta(h))_{X,Y,Z}
=}\\
&=&\phi_{X,Y,Z}\circ (c_{X,Y\tprod Z}\circ(\id\otimes c_{Y,Z}))^{-1}
\circ h_{X\tprod(Y\tprod Z)}\circ
   c_{X,Y\tprod Z}\circ(\id\tprod c_{Y,Z})\\
  &=&(c_{X,Y}^{-1}\otimes \id)\circ c^{-1}_{X\tprod Y,Z}\circ F(\Phi_{X,Y,Z})
     \circ h_{X\tprod(Y\tprod Z)}\circ c_{X,Y\tprod Z}\circ
      (\id\otimes c_{Y,Z})
     \\
\lefteqn{((\Delta\otimes1)\Delta(h)\phi)_{X,Y,Z}=}\\
&=& (c^{-1}_{X,Y}\otimes\id)\circ c^{-1}_{X\tprod Y,Z}\circ
    h_{(X\tprod Y)\tprod Z}\circ
     c_{X\tprod Y,Z}\circ(c_{X,Y}\otimes\id)\circ\phi_{X,Y,Z}\\
   &=&
   (c^{-1}_{X,Y}\otimes1)c^{-1}_{X\tprod Y,Z}h_{(X\tprod Y)\tprod Z}
     F(\Phi_{X,Y,Z})c_{X,Y\tprod Z}(1\otimes c_{Y,Z})\\
\end{eqnarray*}
Both expressions are the same because of naturality:
"$F(\Phi)h=hF(\Phi)$"
This shows quasi coassociativity. For
tensor functors this reduces to coassociativity
and for weak quasi tensor functors  $\phi$ remains quasi invertible.

For the proof of
$(\id\otimes \id\otimes\Delta)(\phi)\cdot(\Delta\otimes \id\otimes \id)(\phi)=
(1\otimes\phi)(\id\otimes\Delta\otimes \id)(\phi)(\phi\otimes1)$.
we refer to Majid's original work \cite{majid2}.

$F$ is a functor between rigid braided tensor categories.
There are isomorphisms $d_X:F(X)^\ast\iso F(X^\ast)$ and
$d_X^\ast:F(X^\ast)^\ast\iso F(X)$. They are used in the definition
of the antipode:
\begin{equation}\label{antidef}
(Sh)_X:=d_X^\ast\circ(h_{X^\ast})^\ast\circ d_X^{\ast-1}
\end{equation}
The proof of the antipode identity will be given in lemma
\ref{antipodenlemma}.

$H$ is quasitriangular by means of  $R\in H\otimes H$:
\begin{equation}	\label{rdef}
R_{X,Y}:=\Psi^{{\Vec}-1}_{F(X),F(Y)}\circ c^{-1}_{Y,X}\circ
 F(\Psi_{X,Y})\circ c_{X,Y}
\end{equation}

$R$ relates the coproduct and the opposite coproduct:
\begin{eqnarray*}
(R\Delta(h)R^{-1})_{X,Y}&=&\Psi^{{\Vec}-1}_{F(X),F(Y)}\circ c^{-1}_{Y,X}\circ
F(\Psi_{X,Y})\circ c_{X,Y}\circ c_{X,Y}^{-1}\circ h_{X\tprod Y}\circ c_{X,Y}\\
&&c_{X,Y}^{-1}\circ F(\Psi_{X,Y})^{-1}\circ c_{Y,X}\circ
\Psi^{{\Vec}}_{F(X),F(Y)}\\
&=& \Psi^{{Vec}-1}_{F(X),F(Y)}\circ c^{-1}_{Y,X}\circ
 F(\Psi_{X,Y})\circ h_{X\tprod Y}\circ\\
&&F(\Psi_{X,Y})^{-1}\circ c_{Y,X}\circ \Psi^{{\Vec}}_{F(X),F(Y)}\\
&=& \Psi^{{\Vec}-1}_{F(X),F(Y)}\circ c^{-1}_{Y,X}\circ h_{Y\tprod X)}
\circ c_{Y,X}\circ \Psi^{{\Vec}}_{F(X),F(Y)}=\Delta'(h)_{X,Y}
\end{eqnarray*}
For the proof of the other two quasitriangularity equations
we refer once more to \cite{majid1} and \cite{haring}.
\end{bew}

%\begin{figure}[ht]
%\setlength{\unitlength}{1mm}
%\begin{picture}(150,60)
%\put(0,40){\morbox{$h_X$}}
%\put(0,30){\bigmorbox{$F(f)$}}
%\put(10,40){$=$}
%\put(20,30){\morbox{$h_Y$}}
%\put(20,40){\bigmorbox{$F(f)$}}
%\put(8,25){{\footnotesize Naturality}}

%\put(35,50){\vertexbox{$c_{X,Y}$}}
%\put(37.55,40){\hugemorbox{$h_{X\tprod Y}$}}
%\put(35,30){\invvertexbox{$c_{X,Y}^{-1}$}}
%\put(37,25){{\footnotesize Coproduct}}

%\put(60,50){\vertexbox{$c_{X,Y}$}}
%\put(60,40){\zopf}
%\put(60,30){\invvertexbox{$c_{X,Y}^{-1}$}}
%\put(60,20){\symzopf}
%\put(62,15){{\footnotesize R matrix}}
%\put(60,40){\dashbox{1.0}(20,10){}}

%\end{picture}
%\caption{\label{hgraph} Weak quasi Hopf algebra structure on $\Nat(F,F)$}.
%The diagrams are in $\Vec$ while portions in dahsed boxes are $F$
%images of morphisms in $\kat$.
%\end{figure}

\begin{lemma} If $F$ is a weak quasi tensor functor then $H$ is a weak quasi
Hopf algebra.
\end{lemma}
\begin{bew}
The additional axioms (the statements already proven remain true!)
are easily verfied using $cc^{-1}=1, c^{-1}c\neq1$:
For (\ref{ff5}) we calculate:
\begin{eqnarray*}
(\id\otimes\id\otimes\epsilon)(\phi)_{X,Y}&=&
(c^{-1}_{X,Y}\otimes\id)\circ c^{-1}_{X\tprod Y,1}\circ
F(\Phi_{X,Y,1})\circ c_{X,Y}\circ (\id\otimes c_{Y,1})=\\
&=&c^{-1}_{X,Y}\circ c_{X,Y}=\Delta(1)_{X,Y}
\end{eqnarray*}
And for (\ref{ff4}):
\begin{eqnarray*}(R^{-1}R)_{X,Y}&=&
c_{X,Y}^{-1}\circ F(\Psi_{X,Y}^{-1})\circ c_{Y,X}\circ
\Psi^{\Vec}\circ\Psi^{\Vec-1}\circ c^{-1}_{Y,X}\circ
F(\Psi_{X,Y})\circ c_{X,Y}\\
&=&c_{X,Y}^{-1}\circ c_{X,Y}=\Delta(1)_{X,Y}
\end{eqnarray*}
Similarly one gets
 (\ref{ff2}):
\begin{eqnarray*}
\phi_{X,Y,Z}\circ \phi^{-1}_{X,Y,Z}
&=&(c_{X,Y}^{-1}\otimes \id)\circ c_{X\otimes Y,Z}^{-1}\circ
F(\Phi_{X,Y,Z})\circ c_{X,Y\otimes Z}\circ (\id\otimes c_{Y,Z})\\&&
(\id\otimes c_{Y,Z}^{-1})\circ c^{-1}_{X,Y\otimes Z}\circ
F(\Phi_{X,Y,Z})^{-1}\circ c_{X\otimes Y,Z}\circ (c_{X,Y}\otimes \id)\\
&=&(c_{X,Y}^{-1}\otimes \id)\circ c_{X\otimes Y,Z}^{-1}\circ
c_{X\otimes Y,Z}\circ (c_{X,Y}\otimes \id)\\
&=&(c_{X,Y}^{-1}\otimes \id)
\Delta(1)_{X\otimes Y,Z}\circ(c_{X,Y}\otimes\id)\\
&=&((\Delta \otimes \id) \Delta(1))_{X,Y,Z}.
\end{eqnarray*}
(\ref{ff1}) is proven in the same way, just as (\ref{ff3}).
\end{bew}

\begin{lemma}
The vector spaces  $F(X)$ are  representation spaces of $H$.
The functor $G:\kat\pfeil\Rep(H)$ is a full tensor functor.
\end{lemma}
\begin{bew}
The representations are
$\varrho_X(h).v:=h_X(v)\quad h\in H,v\in F(X)$.
This induces a functor $G:\kat\pfeil\Rep(H)$.
Morphisms $f\in{\rm Mor}(X,Y)$ are mapped to intertwiners $G(f)=F(f)$:
$G(f)\circ\varrho_X(h)=F(f)\circ h_X=h_Y\circ F(f)=\varrho_Y(h)\circ G(f)$.
$G$ is a tensor functor:
\begin{eqnarray*} (G(X)\tprod G(Y))(h)&=&
(\varrho_X\otimes\varrho_Y)(\Delta(h))=\Delta(h)_{X,Y}=\\
&&c^{-1}_{X,Y}\circ h_{X\tprod Y}\circ c_{X,Y}=
c^{-1}_{X,Y}\circ G(X\tprod Y)\circ c_{X,Y}
\end{eqnarray*}
Here the $c_{X,Y}$ are as maps of vector spaces the same as the $c_{X,Y}$ of
the functor $F$, but because the tensor product on the lefthand side of  this
equation is in the representation category of $H$ they are restricted to the
representation subspace and are therfore isomorphisms. The definitions of
$R$ and $\phi$ are precisley the statements that $G$ is compatible with
associativity and braid isomorphisms.

$G$ is full, because every morphism $T$ in $\Rep(H)$
($T\varrho_Y=\varrho_XT$)
is a constraint that can only exist if it is of the form $T=F(f)$.
\end{bew}

\begin{lemma} \label{lst1}
Let $X,Y\in{\rm Obj}(\kat)$ and $h\in H$.
If $X$ and $Y$ are isomorphic then $h_X$ is determined uniquely by $h_Y$.
If $\kat$ is semisimple then $h\in H$ is determined by its values on $\nabla$
where it may take arbitrary values.
\end{lemma}
\begin{bew}
If $\phi\in{\rm Mor}(X,Y)$ is iso then the naturality condition can be
expressed as $h_Y=F(\phi)\circ h_X\circ F(\phi^{-1})$.

Let $h$ be defined on $\nabla$. Since we assume $\kat$ to be semisimple,
every object is isomorphic to a direct sum of objects in $\nabla$.
By the above remark $h$ is therfore uniquely defined on all objects if it is
uniquely defined on direct sums.
Consider $\bigoplus_{i\in I}X_i, X_i\in\nabla$.
We have morphisms $p_j\in{\rm Mor}(\bigoplus_i X_i,X_j)$ and
$q_j\in{\rm Mor}(X_j,\bigoplus_i X_i)$ such that
$\id_{\bigoplus X_i}=\sum_j q_j\circ p_j$.
Naturality implies
$F(p_j)\circ h_{\bigoplus X_i}=h_{X_j}\circ F(p_j)$ and hence
we have
\[h_{\bigoplus X_i}=F(\sum_j q_j\circ p_j)\circ h_{\bigoplus X_i}
=\sum_j F(q_j)\circ F(p_j)\circ h_{\bigoplus X_i}
=\sum_j F(q_j)\circ h_{X_j} \circ F(p_j)
\]
On different objects in $\nabla$ the function $h$
may take arbitrary values because there
are no morphisms (and hence no naturality constraints)
between inequivalent irreducible objects in an abelian
category.
\end{bew}

\begin{lemma} $G$ is surjective in the sense that it hits
every class of irreps of $H$
\end{lemma}
\begin{bew}
We use lemma \ref{lst1}. It shows that $H$ is a direct sum of
full matrix algebras $M_n(\KK)$.
Each of them has only one irrep. And so $H$ has no other irreducible
representations, because all representations have to reflect commutativity
of the summands and must therefore annihilate all summands but one.
Therefore $H$ has no more irreducible representations classes than $\kat$ has
irreducible object classes.
\end{bew}

\begin{lemma} Faithfulness of $F$ implies that inequivalent objects
yield inequivalent representations.
\end{lemma}
\begin{bew}
Assume $X,Y$ to be inequivalent
objects which are mapped to equivalent representations, i.e.
$F(X)=F(Y),\forall h\in H,
h_X=\varphi\circ h_Y\circ\varphi^{-1}$ with an isomorphism
$\varphi:F(X)\pfeil F(Y)=F(X)$. So the value of
$h$ on $X$ is determined uniquely
by its value on $Y$. This can be done by naturality only if
$\exists f\in{\rm Mor}(X,Y)\exists g\in{\rm Mor}(Y,X)$
such that $F(f)=\varphi,
F(g)=\varphi^{-1}$.
But then (by faithfulness) $f$ and $g$ are iso ($\id_{F(Y)}=F(f)F(g)=F(fg)$;
because of faithfulness only $\id_Y$ is mapped to $\id_{F(Y)}$ and hence
$f=g^{-1}$) contracting our hypothesis.

\end{bew}

\begin{lemma}\label{lllk}
$F(\phi^i\CVO{X}{Y}{Z})\circ c$
form a basis of morphisms in ${\rm Mor}(F(X)\otimes F(Y),F(Z))$
if $F$ is faithful.
\end{lemma}
\begin{bew}
They are linearly independent: Assume $\sum_i\alpha_i
F(\phi^i\CVO{X}{Y}{Z})\circ c$.
By surjectivity of $c$ and linearity of $F$ this implies
$0=F(\sum_i\alpha_i \phi^i\CVO{X}{Y}{Z}))$
and faithfulness of $F$ yields a contradiction.
Further they span the whole space since $G$ is full.
\end{bew}

\begin{lemma} If $F$ is faithful then
$\kat$ and $\Rep(H)$ have the same structural constants
and are therefor equivalent as braided tensor categories.
\end{lemma}
\begin{bew}
Describe $\kat$ as in subsection \ref{peq}.
According to this presentation we have for $X,Y,Z\in\nabla$
matrices $\Omega$ that satisfy $\phi^i\CVO{X}{Y}{Z}\circ\Psi_{Y,X}=\sum_j
\Omega_{i,j}\phi^j\CVO{Y}{X}{Z}$. We apply $F$, multiply $c$ from the right,
introduce $1=cc^{-1}$ and use linearity of $F$ to get
$F(\phi^i\CVO{X}{Y}{Z})\circ c\circ c^{-1}\circ F(\Psi_{Y,X})\circ c=\sum_j
\Omega_{i,j}F(\phi^j\CVO{Y}{X}{Z})\circ c$.
Taking lemma \ref{lllk} into account
and observing that $c^{-1}\circ F(\Psi)\circ c$ is
(by (\ref{zopfhopf}) and (\ref{rdef})) nothing than $\Psi^{\Rep(H)}$
this shows
that $\kat$ and $\Rep(H)$ have the same structure constants.
\end{bew}

\begin{figure}[ht]
\setlength{\unitlength}{1mm}
\begin{picture}(150,105)
\put(2,90){\makebox(8,6){{\footnotesize $F(X)$}}}
\put(12,90){\makebox(8,6){{\footnotesize $F(1)$}}}
\put(0,80){\morbox{$\alpha$}}
\put(10,80){\morbox{$h_1$}}
\put(20,85){$=$}
\put(25,90){\coevpic}
\put(50,90){\line(0,1){10}}
\put(35,80){\evbox{$\ev^\Rep$}}
\put(30,70){\line(0,1){20}}
\put(40,70){\morbox{$h_1$}}
\put(45,67){{\footnotesize $F(1)$}}
\put(30,67){{\footnotesize $F(X)$}}
\put(55,85){$=$}
\put(55,90){\coevpic}
\put(80,80){\line(0,1){20}}
\put(65,80){\morbox{$d_X$}}
\put(65,70){\vertexbox{$c_{X^\ast,X}$}}
\put(70,60){\bigmorbox{$F(\ev)$}}
\put(60,50){\line(0,1){40}}
\put(70,50){\morbox{$h_1$}}
\put(82,85){$=$}
\put(85,90){\coevpic}
\put(110,80){\line(0,1){20}}
\put(100,80){\line(0,1){10}}
\put(85,80){\morbox{$d_X^\ast$}}
\put(95,70){\vertexbox{$c_{X^\ast,X}$}}
\put(100,50){\bigmorbox{$F(\ev)$}}
\put(90,50){\line(0,1){30}}
\put(97.5,60){\hugemorbox{$h_{X^\ast\tprod X}$}}
\put(111,85){$=\sum_i$}
\put(120,90){\coevpic}
\put(145,90){\line(0,1){10}}
\put(140,80){\bigmorbox{$h_i{(2)}$}}
\put(130,80){\morbox{$d_X$}}
\put(120,80){\bigmorbox{$d_X^{\ast-1}$}}
\put(120,70){\bigmorbox{$h_i^{(1)\ast}$}}
\put(120,60){\bigmorbox{$d_X^\ast$}}
\put(130,70){\vertexbox{$c_{X^\ast,X}$}}
\put(135,60){\bigmorbox{$F(\ev)$}}
\put(140,50){\line(0,1){10}}
\put(125,50){\line(0,1){10}}
\put(0,20){$=\sum_i$}
\put(10,35){\coevpic}
\put(20,25){\morbox{$d_X$}}
\put(35,25){\line(0,1){10}}
\put(25,5){\bigmorbox{$F(\ev)$}}
\put(20,15){\vertexbox{$c_{X^\ast,X}$}}
\put(30,35){\bigmorbox{$h_i^{(2)}$}}
\put(15,15){\line(0,1){20}}
\put(7.5,5){\hugemorbox{$S(h_i^{(2)})$}}
\put(40,20){$=\sum_i$}
\put(47.5,5){\hugemorbox{$S(h_i^{(2)})$}}
\put(50,15){\morbox{$\alpha$}}
\put(50,25){\bigmorbox{$h^{(2)}_i$}}
\end{picture}
\caption{\label{antifig}
Proof of the antipode identity as a graphical calculation in $\Vec$}
\end{figure}
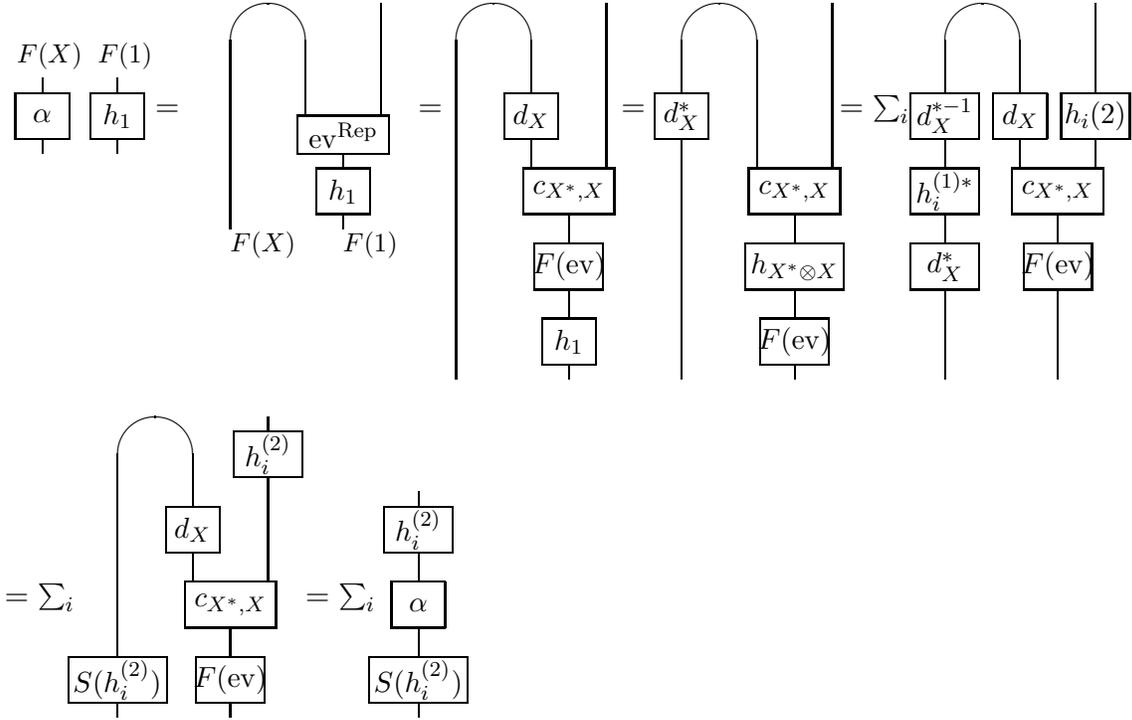

\begin{lemma} \label{antipodenlemma}
$\Rep(H)$ is rigid if $F$ is faithful.
The antipode {\rm (\ref{antidef})} satisfies {\rm (\ref{antiax})}
\end{lemma}
\begin{bew}
$\Rep(H)$ is rigid. The dual representation of $\varrho_X$ is
given by $(\varrho_X)^\ast(h):=(\varrho_X(S(h)))^\ast$ acting on
$F(X)^\ast$.
Note that $\ast$ on the lefthand side of this definition
is the duality in $\Rep$ while on the ridehand side it is the duality in
$\Vec$.

Evaluation and coevaluation are given by
\begin{eqnarray*}
\ev^\Rep_{\varrho_X} &=& F(\ev_X)\circ c_{X^\ast,X}\circ(d_X\otimes \id)\\
\coev^\Rep_{\varrho_X} &=& (\id\otimes d_X^{-1})\circ c^{-1}_{X,X^\ast}\circ
F(\coev_X)
\end{eqnarray*}

We verify the intertwinig property for $\ev^\Rep$ (The proofs for
the coevaluation are identical up to duality symmetry
and are not displayed.):
\begin{eqnarray*}
\lefteqn{\ev^\Rep_{\varrho_X}\circ(\varrho^\ast_X\tprod\varrho_X)(h)=}\\
&&= F(\ev_X)\circ c_{X^\ast,X}\circ (d_X\otimes \id) \circ\\&&
  ((d_X^\ast\otimes \id)\circ (c^{-1}_{X^\ast,X}\circ h_{X^\ast\tprod X}
   \circ c_{X^\ast,X})^{\ast\otimes \id}\circ (d_X^{\ast-1}\otimes \id))^
   {\ast\otimes \id}\\
&&=  F(\ev_X)\circ c_{X^\ast,X}\circ (d_X\otimes \id) \circ
  (d_X^{-1}\otimes \id)\circ c^{-1}_{X^\ast,X}\circ h_{X^\ast\tprod X}
  \circ c_{X^\ast,X}\circ(d_X\otimes \id)\\
&&=F(\ev_X)\circ h_{X^\ast\tprod X}
  \circ c_{X^\ast,X}\circ(d_X\otimes \id)\\
&&=h_1\circ F(\ev_X)\circ c_{X^\ast,X}\circ(d_X\otimes \id)
= \varrho_1(h)\circ \ev^\Rep_{\varrho_X}
\end{eqnarray*}
Because $\Rep$ is in general not strict (even if $\kat$ is strict
($\Phi^\kat =\id$) which we will assume) we have to insert an associator
into the  the fundamental $\ev/\coev$ property (\ref{evcoevax}):
\begin{eqnarray*}
\lefteqn{(\id\otimes \ev^\Rep_{\varrho_X})\circ \Phi^{\Rep-1}\circ
(\coev^\Rep_{\varrho_X}\otimes \id)=}\\
&&= (\id\otimes \ev^\Rep_{\varrho_X})\circ
  (\varrho_X\tprod\varrho_X^\ast\tprod\varrho_X)(\phi^{-1})\circ
  (\coev^\Rep_{\varrho_X}\otimes \id)\\
&&=(\id\otimes F(\ev_X))\circ (\id\otimes c_{X^\ast,X})\circ
     (\id\otimes d_X\otimes \id)\circ\\
&&   (\id\otimes d_X^{-1}\otimes \id)\circ(\id\otimes c^{-1}_{X^\ast,X})\circ
    c^{-1}_{X,X^\ast\tprod X}\circ F(\Phi^{\kat-1})\circ
  c_{X\tprod X^\ast,X}
  \circ (c_{X,X^\ast}\otimes \id)\circ \\
&& (\id\otimes d_X\otimes \id)\circ
  (\id\otimes d^{-1}_X\otimes \id)\circ (c^{-1}_{X,X^\ast}\otimes \id)\circ
  (F(\coev_X)\otimes \id)\\
&&= (\id\otimes F(\ev_X))\circ (c^{-1}_{X,X^\ast\tprod X})\circ
  c_{X\tprod X^\ast,X}\circ (F(\coev_X)\otimes \id)\\
&&=c^{-1}_{X,1}\circ(\id\otimes F(\ev_X))\circ(F(\coev_X)\otimes \id)
\circ c_{1,X}
=\id
\end{eqnarray*}
The antipode identities involve elements $\alpha,\beta\in H(\kat,F)$.
\begin{eqnarray}
\alpha_X& := & (\id\otimes \ev^\Rep_{\varrho_X})\circ(\coev^\Vec_{F(X)}
   \otimes \id):F(X)\pfeil F(X)\\
\beta_X& := & (\id\otimes \ev^\Vec_{F(X)})\circ(\coev^\Rep_{\varrho_X}
   \otimes \id):F(X)\pfeil F(X)
\end{eqnarray}
This gives well defined elements in $H$:
$F(f)\circ\alpha_X=\alpha_Y\circ F(f)$ holds because $d, c$ are functorial.

Applying $\ev^\Vec\otimes id$ yields
$\ev_{\varrho_X}^\Rep=\ev^\Vec\circ(\id\otimes\alpha_X)$.
Obviously $\alpha=\beta^{-1}$ if $\Rep(H)$ is strict.
The proof of the antipode identity is given in the following calculation
and additionally in graphical notation in figure \ref{antifig}.
\begin{eqnarray*}
\lefteqn{(\epsilon(h)\alpha)_X=}\\
&&=(\id\otimes h_1)\circ(\id\otimes F(\ev_X))\circ
   (\id\otimes c_{X^\ast,X})\circ
   (\id\otimes d_X\otimes \id)\circ(\coev^\Vec_{F(X)}\otimes \id)\\
&&=(\id\otimes F(\ev_X))\circ(\id\otimes h_{X^\ast\tprod X})\circ
   (\id\otimes c_{X^\ast,X})\otimes(d_X^\ast\otimes \id\otimes \id)\circ
   (\coev^\Vec\otimes \id)\\
&&=(d_X^\ast\otimes \id\otimes \id)\circ(\id\otimes
    (F(\ev_X)\circ c_{X^\ast,X}))
    \circ(\id\otimes (d_X\circ d_X^{-1})\otimes \id)\circ\\
&&  \quad\sum_i (\id\otimes h_i^{(1)}\otimes h_i^{(2)})
\circ(\coev^\Vec\otimes \id)\\
&&=\sum_i d_X^\ast\circ h_i^{(1)\ast}\circ d_X^{-1\ast}\circ
   (\id\otimes F(\ev_X)\circ c_{X^\ast,X})\circ(\id\otimes d_X\otimes\id)\circ
    (\coev^\Vec\otimes h_i^{(2)})\\
&&=\left( \sum_i S(h_i^{(1)})\alpha h_i^{(2)}\right)_X
\end{eqnarray*}
The second antipode axiom (\ref{antiax2}) involving $\beta$ is
established similarly.
\end{bew}

\begin{lemma} If $F$ is isometric then $H$ is involutive
and the representations are unitary:
$\varrho(h^\dagger)=\varrho(h)^\dagger$.
\end{lemma}
\begin{bew}
The involution is given by: $(h^\dagger)_X:=(h_X)^\dagger$.
Applying $F\circ\dagger=\dagger\circ F$ to the naturality condition
implies that $H$ is closed under this operation.
Multiplicativity carries over from vector space endomorphisms.
$\Delta(h^\dagger)=\Delta(h)^\dagger$ follows easily from the fact
that $c^\dagger=c^{-1}$:
$\Delta(h)^\dagger_{X,Y}=(c_{X,Y}^{-1}\circ
h_{X\tprod Y}\circ c_{X,Y})^\dagger=
c_{X,Y}^{-1}\circ h_{X\tprod Y}^\dagger \circ c_{X,Y}$.

The proof of the compatibility of the involution $\dagger$ and the antipode
$S$ uses the fact that in the category of finite dimensional Hilbert
spaces the duality map is given by $V^\ast=V$ and $\ev/\coev$ is given by
the scalar product. Using this one calculates
\[(S(h)^\dagger)_X=(S(h)_X)^\dagger=(d_X^\ast\circ h^\ast_{X^\ast}\circ
 d_X^{\ast-1})^\dagger=d_X^{\ast-1\dagger}\circ h_{X^\ast}^{\ast\dagger}\circ
 d_X^{\ast\dagger}=d_X^\ast\circ h_{X^\ast}^{\dagger\ast}\circ d_X^{\ast-1}
=S(h^\dagger)\]
\end{bew}

\begin{lemma} If $\kat$ is a ribbon category then $H$ is a ribbon
Hopf algebra in the sense of \cite{resh}.
\end{lemma}
\begin{bew}
The ribbon element $v\in H$ is defined by $v_X:=F(\sigma(X))$. It is
central because $\sigma$ is functorial. Further one calculates:
$\epsilon(v)=F(\sigma(1))=1$ and
\begin{eqnarray*}
S(v)_X=d^\ast_X\circ (v_{X^\ast})^\ast\circ d_X^{\ast-1}&=&
d^\ast_X\circ F(\sigma(X^\ast))^\ast\circ d_X^{\ast-1}=\\&&
d^\ast_X\circ F(\sigma(X)^\ast)^\ast\circ d_X^{\ast-1}=
F(\sigma(X))=v_X
\end{eqnarray*}
\begin{eqnarray*}
\Delta(v)_{X,Y}&=&c^{-1}_{X,Y}\circ
 F(\sigma(X\tprod Y))\circ c_{X,Y}\\
 &=&c^{-1}_{X,Y}\circ F(\Psi_{X,Y}^{-1})\circ c_{Y,X}\circ c_{Y,X}^{-1}\circ
F(\Psi_{Y,X}^{-1})\circ c_{X,Y}\circ c_{X,Y}^{-1}\circ \\&&
F(\sigma(X)\tprod \sigma(Y))\circ c_{X,Y}\\
&=&((R_{2,1}R)^{-1})_{X,Y}\circ (F(\sigma(X))\otimes
F(\sigma(Y)))\\
&=&((R_{2,1}R)^{-1}(v\otimes v))_{X,Y}
\end{eqnarray*}
\end{bew}

\begin{lemma}\label{ultralemma}
If $F$ is an ultra weak quasi tensor functor then
$H$ is an ${\rm End}(F(1))$ ultra weak quasi Hopf algebra.
\end{lemma}
\begin{bew}
The bimodule actions are defined to be:
\begin{equation}
\mu_l(a\otimes h)_X:=c_{1,X}\circ(a\otimes h)\circ c^{-1}_{1,X}\quad
 \mu_r(h\otimes a)_X:=c_{X,1}\circ(h\otimes a)\circ c^{-1}_{X,1}\quad
 a\in{\rm End}(F(1))
 \end{equation}
The definition of $\epsilon$ doesn't have to be changed but
the unit is now defined more general to be
\begin{equation}
\eta(a):=\mu_l(a\otimes 1)=\mu_r(1\otimes a)
\end{equation}
The counit property is fulfilled:
\[(\mu_l(\epsilon\otimes \id)\Delta(h))_X=c_{1,X}\circ c_{1,X}^{-1}\circ
h_{1\tprod X}\circ c_{1,X}\circ c_{1,X}^{-1}=h_X\]
\end{bew}

Collecting results together we have:

\begin{theo}[Generalized Majid's reconstruction theorem]
\label{majidrt1}
Let $\kat$ be a rigid braided tensor category and
$F:\kat\pfeil\Vec$
a weak quasi tensor functor. Then the set $H=\Nat(F,F)$ carries the structure
of a weak quasi Hopf algebra and there is a
functor $G:\kat\pfeil \Rep(H)$ such that
$\kat\stackrel{G}{\pfeil}\Rep(H)\stackrel{V}{\pfeil}\Vec$
composes to $F$.
$G$ maps inequivalent objects to inequivalent representations
if $F$ is faithful.
$G$ is full. $G$ is faithful iff $F$ is faithful.
Hence in the case of a faithful functor and a semisimple category,
$\kat$ and $\Rep(H)$ are equivalent braided tensor categories.
$\Rep(H)$ is rigid if $F$ is faithful and it is $C^\ast$ if $\kat$ is so.
The structure matrices {\rm (}see section {\rm \ref{peq}}{\rm)}
of $\kat$ and $\Rep(H)$ coincide.
A ribbon structure on $\kat$ induces a ribbon structure on $H$ \cite{resh}.
The structure of $H$ is determined by $F$:

$F$ is tensor functor \quad$\Longrightarrow$\quad
$H$ is quasitriangular Hopf algebra

$F$ is quasi tensor functor \quad$\Longrightarrow$\quad
$H$ is quasitriangular quasi Hopf algebra

$F$ is weak quasi tensor functor \quad$\Longrightarrow$\quad
$H$ is quasitr. weak quasi Hopf algebra

$F$ is ultra weak quasi tensor functor \quad$\Longrightarrow$\quad
$H$ is quasitr. ultra weak quasi Hopf algebra
\end{theo}

Combinig this with the construction of weak quasi tensor functors we
conclude:

\begin{kor} Every rational semisimple rigid braided tensor category
is the representation category of some weak quasi Hopf algebra.
\end{kor}

\subsection{Questions of non uniqueness}

The reconstruction of $H$ from a given category $\kat$
presented in this paper is not unique. It can be checked
that in some typical examples  there is an infinite number of
weak dimension functions. But there is even more freedom
because of the choice of epimorphisms $C$ in the proof of
proposition \ref{funktorkon}.

\begin{bem}
Let $F,\schlange{F}:\kat\pfeil{\Vec}$
denote two faithful (weak) quasi  tensor functors
constructed as in proposition {\rm \ref{funktorkon}}
by the same dimension function.
Then the reconstructed (weak) quasi Hopf algebras $H$ and $\schlange{H}$
are equal up to twist equivalence (in the sense of Drinfeld).
\end{bem}
\begin{bew}
$F$ and $\schlange{F}$ differ only by differnt choices of $C$.
However because they share the same dimension function there
is a family of isomorphisms $\varphi$ such that
$\schlange{c}_{X,Y}=\varphi_{X,Y}\circ c_{X,Y}$.
$H$ and $\schlange{H}$ are then equal as algebras.
Their coalgebra structure however differs.
\begin{eqnarray*}
\schlange{\Delta}(h)_{X,Y}&=&\schlange{c}_{X,Y}^{-1}\circ
h_{X\tprod Y}\circ \schlange{c}_{X,Y}\\
&=&c^{-1}_{X,Y}\circ \varphi^{-1}_{X,Y}\circ c_{X,Y}\circ c_{X,Y}^{-1}\circ
h_{X\tprod Y}\circ c_{X,Y}\circ c_{X,Y}^{-1}\circ \varphi_{X,Y}\circ  c_{X,Y}\\
&=& T_{X,Y}\circ \Delta(h)_{X,Y}\circ T_{X,Y}^{-1}\\
&=&(T\Delta(h)T^{-1})_{X,Y}
\end{eqnarray*}
Here we have inserted $1=cc^{-1}$ twice and introduced
the twist element $T\in H\otimes H$, defined by
$T_{X,Y}:=c^{-1}_{X,Y}\circ\varphi_{X,Y}^{-1}\circ c_{X,Y}$.
(In the weak case $T$ is not invertible,
but one has $TT^{-1}=T^{-1}T=\Delta(1)$.)
Note that $T$ really is an element of $H\otimes H$ because the
dependence of $T_{X,Y}$ on $X,Y$  obeys
the naturality condition as is easily seen from the
definition of $F$ and $c$.

Of course the $R$ element gets twisted alike:
\begin{eqnarray*}
\schlange{R}_{X,Y}
&=& \Psi^{\Vec-1}_{F(Y),F(X)}\circ\schlange{c}_{Y,X}^{-1}
\circ F(\Psi_{X,Y})\circ\schlange{c}_{X,Y} \\
&=&\Psi^{\Vec-1}_{F(Y),F(X)}\circ c_{Y,X}^{-1}\circ
\varphi_{Y,X}^{-1}\circ F(\Psi_{X,Y})\circ
\varphi_{X,Y}\circ c_{X,Y}\\
&=& \Psi^{\Vec-1}_{F(Y),F(X)}\circ c_{Y,X}^{-1}\circ\varphi_{Y,X}^{-1}
\circ c_{Y,X}\circ \Psi^{\Vec}_{F(X),F(Y)}\circ\Psi^{\Vec-1}_{F(X),F(Y)}
\circ c_{Y,X}^{-1}\circ\\
&& F(\Psi_{X,Y})\circ c_{X,Y}\circ c_{X,Y}^{-1}\circ
\varphi_{X,Y}\circ c_{X,Y}\\
&=&\Psi^{\Vec-1}_{F(X),F(Y)}\circ T_{Y,X}\circ\Psi^{\Vec}_{F(X),F(Y)}
\circ R_{X,Y}\circ t^{-1}_{X,Y}\\
&=&(T_{2,1}RT_{1,2}^{-1})_{X,Y}
\end{eqnarray*}
A similar calculation shows that
\[\schlange{\phi}=T_{1,2}^{-1}(\Delta\otimes \id)(T^{-1})
\phi(\id\otimes\Delta)(T)T_{2,3}\]
\end{bew}

The results of the previous remark can be nicely interpreted in the language
of nonabelian cohomology \cite{majidbuch} where the n-cochains are given
by the invertible elements in $H^{\otimes n}$ and the coboundary operator is
defined to be
\[\delta(\gamma):=\prod_{i=1,3,...}\Delta_i(\gamma)\prod_{i=0,2,..}
\Delta_i(\gamma)^{-1}\quad\in H^{\otimes n+1}\qquad\gamma\in H^{\otimes n}\]
with $\Delta_0(\gamma):=1\otimes\gamma,\Delta_{n+1}(\gamma):=\gamma\otimes1,
\Delta_i(\gamma):=(\id\otimes...\otimes\Delta\otimes...\otimes \id)(\gamma)$.
Then the pentagon identity for $\phi\in H^{\otimes 3}$ is the statement
that $\phi$ is a 3-cocycle, $\delta(\phi)=1\otimes1\otimes1\otimes1$.
Returning to the previous remark we see that $\schlange{\phi}$ can be made
trivial iff $\phi$ is a coboundary $\phi=\delta(T)$.

{\em\bf Acknowledgement:} It is a pleasure to thank all those who
helped in the trail of this work, especially
F. Constantinescu, T. tom Dieck, H. Kratz, G. Mack, S. Majid,
J. E. Roberts and A. Schmidt.

Thanks to the Stu\-dien\-stift\-ung des deutschen Volkes and the Deu\-tsche
For\-schungs\-gemeinschaft, SFB 170,  for
financial support and more.

\twocolumn


\small
\begin{thebibliography}{99}
  \bibitem{deligne} P. Deligne, J. S. Milne: Tannakian Categories,
  in P. Deligne et al,
  Springer LNM 900
  \bibitem{drinfeld} V.G. Drinfeld:
  Quasi-Hopf Algebras and Knizhnik-Zamolodchikov Equations,
  in A. A. Belavin: Problems of Modern QFT, Springer 1989
  \bibitem{Haag} R. Haag, Local Quantum Physics, Springer 1992
  \bibitem{haring} R. H\"aring:
  Quanten-Symmetrie, diploma thesis, Frankfurt 1993
  \bibitem{kw} D. Kazhdan, H. Wenzl: Reconstruction of Monoidal Categories,
  Adv. in Soviet. Math. 16, 11-136 (1993)
  \bibitem{kerler} T. Kerler:
  Non-Tannakian Categories in QFT, Preprint ETH-TH/91-51
  \bibitem{kratz} H. Kratz: Weak quantum symmetry in all minimal models,
   Physics Letters B
  \bibitem{macl} S. MacLane:
  Categories for the working Mathematician, Springer LNM
  \bibitem{ms0} G. Mack, V. Schomerus:
  Conformal Field Algebras with Quantum Symmetry
  from the Theory of Superselection Sectors,
  Comm. Math. Ph. 134, 139-196 (1990)
  \bibitem{ms1} G. Mack, V. Schomerus: Quasi Hopf quantum symmetry in
                quantum theory, Nucl. Ph. B370, 185-230 (1992)
  \bibitem{majidbuch} S. Majid: Foundations of Quantum Group Theory,
     in preparation
  \bibitem{majid} S. Majid:
  Quaitriangular Hopf Algebras and Yang-Baxter-Equations,
    Int. J. Mod. Ph. A5(1), 1-91 (1990)
  \bibitem{majid1} S. Majid:
  Reconstruction Theorems and RCQFT, Int. J. of Mod. Ph. A6(24),
  4359-4374 (1991)
  \bibitem{majid2} S. Majid:
  Quasi-Quantum Groups as Internal Symmetries of Topological
  QFTs, Lett. Math. Ph. 22, 83-90 (1991)
  \bibitem{majid3} S. Majid: Braided Groups and
  Algebraic QFTs, Lett. Math. Ph. 22, 167-176 (1991)
  \bibitem{majid10} S. Majid: Tannaka-Krein Theorems for Quasi-Hopf Algebras
   and other results, Contemp. Math. 134, 219 (1992)
  \bibitem{moores2} G. Moore, N. Seiberg:
   Classical and Quantum CFT, Comm. Math. Ph. 123, 177-254 (1989)
  \bibitem{nill} F. Nill: Fusion Structure from Quantum Groups II: Why
    Truncation is Necessary, Lett. M. Ph. 29, 83-90
  \bibitem{resh} N. Y. Reshetikhin, V. G. Turaev: Ribbon Graphs and Their
   Invariants derived from Quantum Groups, Comm. Math. Ph. 127, 1-26 (1990)
  \bibitem{schomerus1} V. Schomerus: Quantum-Symmetry in Quantum Theory,
   PhD thesis, Hamburg 1993
  \bibitem{turaev} V. G. Turaev: Quantum Invariants of Knots and 3-Manifolds,
   de Gruyter, Berlin 1994
  \bibitem{ulbrich} K.-H. Ulbrich: On Hopf-Algebras and Rigid
     Monoidal Categories, Israel J. Math. 72, 252-256 (1990)
\end{thebibliography}
\end{document}